\documentclass[a4paper,11pt]{article}
\UseRawInputEncoding
\usepackage{amsmath}
\usepackage{amssymb}
\usepackage{slashbox}
\usepackage{texdraw}
\usepackage{mathrsfs}
\usepackage{epsfig}
\usepackage{amsfonts}
\usepackage{cite}

\usepackage[dvips]{psfrag}
\usepackage{subfigure}
\usepackage{color}
\usepackage{booktabs}

\usepackage{geometry}
\usepackage{hyperref}
\hypersetup{colorlinks=true,
            linkcolor=blue,
            anchorcolor=blue,
            citecolor=blue}

\usepackage{graphicx}

\newtheorem{Definition}{Definition}[section]
\newtheorem{Proposition}{Proposition}[section]
\newtheorem{Lemma}{Lemma}[section]
\newtheorem{Theorem}{Theorem}[section]

\newcommand{\bpf}{{\bf Proof:\ \ }}
\newcommand{\epf}{\mbox{}\hfill $\Box$}
\newcommand{\di}{\displaystyle}

\usepackage{caption}
\usepackage {indentfirst}

\begin{document}
\setlength{\baselineskip}{15.2pt} \pagestyle{myheadings}

\title{Maximum number of limit cycles for Abel equation having coefficients with linear trigonometric functions}

\author{\small Xiangqin Yu$^{a}$, Jianfeng Huang$^{b}$, Changjian Liu{$^{a,*}$} \\
{\small $^{a}$ School of Mathematics (Zhuhai), Sun Yat-sen University,
  Zhuhai, 519082,  P.R. China }\\
{\small $^{b}$ Department of Mathematics, Jinan University, Guangzhou, 510632, P.R. China }}

\renewcommand{\thefootnote}{}
\footnotetext{$^*$Corresponding author.
\par
E-mail addresses: yuxq25@mail2.sysu.edu.cn(X. Yu), thuangjf@jnu.edu.cn(J. Huang),
\par
liuchangj@mail.sysu.edu.cn(C. Liu).}

\date{}
\maketitle
\begin{abstract} This paper devotes to the study of the classical Abel equation $\frac{dx}{dt}=g(t)x^{3}+f(t)x^{2}$, where $g(t)$ and $f(t)$ are trigonometric polynomials of degree $m\geq1$. We are interested in the problem that whether there is a uniform upper bound for the number of limit cycles of the equation with respect to $m$, which is known as the famous Smale-Pugh problem. In this work we generalize an idea from the recent paper (Yu, Chen and Liu, arXiv:$2304.13528$, $2023$) and give a new criterion to estimate the maximum multiplicity of limit cycles of the above Abel equations. By virtue of this criterion and the previous results given by {\'A}lvarez et al. and Bravo et al., we completely solve the simplest case of the Smale-Pugh problem, i.e., the case when $g(t)$ and $f(t)$ are linear trigonometric, and obtain that the maximum number of limit cycles, is three.
\end{abstract}


\vspace{1mm}
\section{Introduction and main result\label{introduction}}
The theories on the Abel differential equations provide tools to study the
problems in many areas of mathematics, especially in the qualitative theory of differential equations. It initially arose in the studies of Abel \cite{abel1829} on the theory of elliptic functions, and appeared in the reduction of order for the high-order differential equations, therefore were frequently found in the varied real models, such as the Li\'{e}nard equations \cite{gine2011}, the evolution equations from the cosmological models \cite{yurov2014application}, and the reaction-diffusion equation that describes the evolution of glioblastomas \cite{harko2015bio}.
In the past decades, the Abel equations have gained increasing attentions because they not only are applied to the characterizations of real periodic phenomena (see e.g. the tracking control problem \cite{fossas2008iterative}, the perturbed pendulums \cite{gasull2016number,gasull2020chebyshev} and the seasonal prey harvesting \cite{xu2005harvesting}), but also play an important role in the study of the Hilbert's $16$th problem (see e.g. \cite{cherkas1976number,li2003hilbert}). The second part of the Hilbert's $16$th problem essentially asks that whether there exists a finite number $\mathscr H(m)$ such that each planar polynomial differential system of degree $\leq m$ has no more than $\mathscr H(m)$ limit cycles.

An Abel equation of generalized type is a one-dimensional non-autonomous differential equation that is written as
\begin{equation}\label{eq1}
  \di\frac{dx}{dt}=\displaystyle\sum_{i=0}^{n} A_{i}(t)x^{i},
\end{equation}
where $x\in \mathbb R,t\in \mathbb R$ and $A_{i}:\mathbb R\rightarrow \mathbb R,i=1,...,n,$ are analytic $2\pi$-periodic functions.
For the cases that $n=3, 2$ and $1$, the equation is called the (classical) Abel equation, the Riccati equation and the linear equation, respectively.
 Let $x(t,x_{0})$ be the solution of equation \eqref{eq1} satisfying $x(0,x_{0})=x_{0}$. We say that $x(t,x_{0})$ is a periodic solution of the equation if $x(2\pi,x_{0})=x_{0}$. And an orbit $x=x(t,x_{0})$ is called a periodic orbit (resp. limit cycle) of the equation, if $x(t,x_{0})$ is a periodic solution (resp. isolated periodic solution).

The studies of the Hilbert's $16$th problem via the equation \eqref{eq1} mainly goes back to the 70's of the last century.
It is shown that there is a series of planar polynomial differential systems that can be reduced to the Abel equations, using e.g. the polar coordinates and Cherkas¡¯
transformation  \cite{devlin1998cubic}. For more details we refer to the works on the quadratic systems (see e.g. \cite{lins1980number,gasull1990limit}), the rigid
systems, and some other general planar polynomial differential systems (see e.g. \cite{gasull1990limit,huang2017estimate}).
For this reason, estimating the number of limit cycles of the Abel equation \eqref{eq1}, becomes a significant problem and is extensively studied by the researchers.

 The first progress in this research area is motivated by Lins-Neto \cite{lins1980number} and Lloyd \cite{lloyd1973number,lloyd1975class,lloyd1979note}. They proved that there exists at most $n$ limit cycle(s) of the equation \eqref{eq1} for $n=1, 2$ (i.e., for the linear type and the Riccati type). However, the situation becomes extremely complicated for $n=3$ (i.e., for the classical Abel type). In \cite{lins1980number} Lins-Neto shown that the maximum number of limit cycles of the classical Abel equations, is unbounded without additional conditions (see also \cite{panov1999diversity}). More concretely, by using the bifurcation methods he presented that the equation of the form
 \begin{align}\label{eq2}
 \frac{dx}{dt}=A_3(t)x^3+A_2(t)x^2
 \end{align}
 with $A_3(t)$ and $A_2(t)$ being the trigonometric polynomials of degree $m$, can possesses at least $m$ limit cycles.
 Such result was generalized later to the case $n>3$ in the work \cite{gasull2006limit}.

 Due to these facts and the above backgrounds, a more specific related version of the Hilbert's $16$th problem arises (see for instance \cite{lins1980number,ilyashenko2000hilbert}): {\em Whether the maximum number of limit cycles for the equation \eqref{eq1} with the coefficients $A_i$'s being trigonometric polynomials of degree $m$, is bounded in terms of $m$}? For the initial case, i.e., the equation of the form \eqref{eq2}, it is known as the Smale-Pugh problem \cite{smale2000mathematical}.

 In fact, this is a very difficult problem that is still open for even the linear trigonometric classical Abel equations.
 So far, there are two lines of nice works, exploring the problem from different perspectives.
 The first line of works is based on the definite sign(s) hypothesis for some coefficients or their linear combinations of the equation, which essentially come from the ¡°transversality¡± of some curve(s) and the orbits of the equation. To illustrate, Lloyd \cite{lloyd1979note} and Pliss \cite{pliss1966nonlocal} (resp. Gasull and Llibre \cite{gasull1990limit}) obtained that the equation \eqref{eq1} with $n=3$ has at most three limit cycles if $A_{3}(t)\neq0$ (resp. $A_2(t)\neq0$ and $A_0(t)\equiv0$). Later Gasull and Guillamon \cite{gasull2006limit} extended such result, showing that the number of limit cycles for the equation $dx/dt=A_{n_1}(t)x^{n_1}+A_{n_2}(t)x^{n_2}+A_{1}(t)x$ with $n_1>n_2>1$ is bounded if either $A_{n_1}(t)$ or $A_{n_2}(t)$ has definite sign. In \cite{alvarez2007new} {\'A}lvarez et al. proved for the first time that the equation \eqref{eq2} can only possess at most one non-zero limit cycle if there exists a linear combination of $A_3(t)$ and $A_2(t)$ with fixed sign. This yields one of the key results that will solve the linear trigonometric case of the smale-pugh problem, which is presented in Theorem \ref{AT1} below.
 Huang and Liang \cite{huang2020geometric} developed the idea in \cite{alvarez2007new} and provided a criterion for the linear trigonometric generalized Abel equations having at most $n$ limit cycles.
 For more works of this line, see \cite{alvarez2015limit,ilyashenko2000hilbert,huang2012periodic,huang2017estimate}. The second line of works mainly focuses on the hypothesis for the symmetries of the coefficients of the equation \eqref{eq1}, see for instance \cite{alvarez2013existence,bravo2009limit,bravo2008nonexistence} and the references therein. It is notable that under such symmetry hypothesis the authors in \cite{bravo2009limit} provide the estimate for the number of limit cycles of the linear trigonometric classical Abel equations, which will be the second key result applied in our paper. We state it in Theorem \ref{AT2}.

For the other related works, the readers are referred to papers \cite{alvarez2017centers,huang2020number}.

Let us come back to the initial case of the problem, i.e., the Smale-Pugh problem. In this paper, we focus on the simplest unsolved case
\begin{equation}\label{eq3}
\frac{dx}{dt}=(a_0+a_1\sin t+a_2\cos t)x^3+(b_0+b_1\sin t+b_2\cos t)x^2,
\end{equation}
where $a_0, a_1, a_2, b_0, b_1, b_2\in\mathbb R$.
The equation \eqref{eq3} was first systematically studied in \cite{alvarez2007new} and \cite{bravo2009limit}. Following the different ideas (as introduced above) the authors of these two works gave several important criteria to estimate the maximum number of limit cycles of the equation, which are summarized as below:
 \begin{Theorem}[\cite{alvarez2007new}]\label{AT1}
Suppose that the equation \eqref{eq3} has not a center at $x=0$.
If there exists $\lambda\in\mathbb R$ such that $\lambda(a_0+a_1\sin t+a_2\cos t)+(b_0+b_1\sin t+b_2\cos t)$ does not vanish identically and does not change sign, i.e., one of the conditions $a_{0}^{2}\geq a_{1}^{2}+a_{2}^{2}$, $b_{0}^{2}\geq b_{1}^{2}+b_{2}^{2},$ or $(a_{1}b_{0}-a_{0}b_{1})^{2}+(a_{0}b_{2}-a_{2}b_{0})^{2}\geq (a_{1}b_{2}-a_{2}b_{1})^{2}$ is satisfied, then it has at most one non-zero periodic orbit. Furthermore, when this periodic orbit exists, it is hyperbolic.
\end{Theorem}

 \begin{Theorem}[\cite{bravo2009limit}]\label{AT2}
If $a_{0}b_{0}=0$, then equation \eqref{eq3} has at most three limit cycles, including $x=0$. Moreover, the upper bound is sharp and in case of having three limit cycles, one is in the region $x>0$ and the other is in the region $x<0$.
\end{Theorem}
We remark that the result in \cite{bravo2009limit} is improved recently by Bravo, Fern{\'a}ndez and Ojeda \cite{bravo2023stability}. The authors show that the maximum number of limit cycles of the equation \eqref{eq3} is still three when $a_{0}b_{0}$ is sufficiently small.
Another new work is given by Yu, Chen and Liu \cite{yu2023number}, in which the equation \eqref{eq3} with $a_2=b_1=0$ reduced from the Josephson equations is studied. The authors also prove that the number of limit cycles of such equation, does not exceed three.

On the other hand, by means of the Hopf (resp. Poincar\'{e}) bifurcation it was obtained in \cite{alvarez2007new} (resp. \cite{huang2020number}) that the equation \eqref{eq3} can has at least three limit cycles. Throughout these works all the evidences yield the following specific version of the simplest case of the Smale-Pugh problem, which is proposed by Gasull and stated as the $6$th of the $33$ open problems in  \cite{gasull2021some}:
\vskip0.3cm
\noindent{\bf Problem (\cite{gasull2021some})} {\em Whether the maximum number of limit cycles of the equation \eqref{eq3}, is three}?
\vskip0.3cm

The purpose of this paper is to study the cases of the equation \eqref{eq3} which are not considered in the previous works \cite{alvarez2007new,bravo2009limit}, and then completely solve the above Smale-Pugh-Gasull problem in \cite{gasull2021some}. By generalizing the recent work in \cite{yu2023number}, we give a criterion for the maximum multiplicity of the non-zero limit cycles of the Abel equation \eqref{eq2} with general coefficients (see section \ref{mA1} for details). Then, applying this criterion, and taking Theorem \ref{AT1} and Theorem \ref{AT2} into account, we finally provide a positive answer as below:
\begin{Theorem}\label{T2}
The equation \eqref{eq3}, with arbitrary parameters $a_0,a_1,a_2,b_0,b_1,b_2\in\mathbb R$, has at most three limit cycles (including $x=0$). Moreover, this upper bound is sharp.
\end{Theorem}

The layout of the rest of this paper is as follows: In Section \ref{pre}, we present some preliminary results. In Section \ref{mA1}, we give a new criterion to estimate the maximum multiplicity of limit cycles of the Abel equation \eqref{eq2}. Finally, the proof of Theorem \ref{T2} is given in Section \ref{mA2}.

\section{Preliminaries \label{pre}}
We recall some basic tools in this section. All of them will be necessary in our argument. The first one characterizes the derivatives of the Poincar\'{e} map for the one-dimensional differential equations, which is initially provided in \cite{lloyd1979note}.

\begin{Proposition}[\cite{lloyd1979note}]\label{pro11}
Suppose that $x(t,x_{0})$ is the solution of the one-dimensional differential equation ${\frac{dx}{dt}}=S(x,t)$ that satisfies $x(t_0,x_{0})=x_{0}$. Denote by $L(x_0)=x(t_0+2\pi, x_0)$. Then
 \begin{enumerate}
\item[\rm{(i)}]  $L^{'}(x_0)=\exp\bigg[\displaystyle{\int_{t_0}^{t_0+2\pi} \frac{\partial S}{\partial x}(x(t, x_0),t)\, dt}\bigg];$
\item[\rm{(ii)}] $L^{''}(x_0)=L^{'}(x_0)\bigg[\displaystyle{\int_{t_0}^{t_0+2\pi} \frac{\partial^{2} S}{\partial x^{2}}(x(t, x_0),t)}\cdot \exp{\displaystyle{\bigg\{\int_{t_0}^{t} \frac{\partial S}{\partial x}(x(\tau, x_0),\tau)} \, d\tau \bigg\}}\, dt\bigg].$
\end{enumerate}
\end{Proposition}

The followings present the concept and several properties of rotated one-dimensional differential equations, which are directly adapted from the classical theory of rotated vector fields (see for instance \cite{duff1953number}). For more details of these adaptions we also refer readers to the works \cite{bravo2015stability,bravo2023stability,han2018theory,yu2023number}.
\begin{Definition}\label{def1}
Consider a family of equations
\begin{equation}\label{eq4}
  \frac{dx}{dt}=S(x,t;\lambda),
\end{equation}
where $t\in E\subset\mathbb R$, $x\in J\subset\mathbb R$ and the parameter $\lambda\in \mathbb R$. Then we say that \eqref{eq4} defines
a family of rotated equations in $E\times J$, if $\frac{\partial{S}}{\partial{\lambda}}(x,t;\lambda)>0$ for $\lambda\in\mathbb R$.
\end{Definition}

\begin{Proposition}\label{pro12}
Suppose that the equation \eqref{eq4} defines a family of rotated equations in the region $[0,2\pi]\times J\subset[0,2\pi]\times\mathbb R$.
\begin{enumerate}
\item[\rm{(i)}] If $x=x(t)$ is a stable or an unstable limit cycle of the equation \eqref{eq4}$|_{\lambda=\lambda_0}$ in the region, then there exists $\delta>0$ such that for $\lambda\in(\lambda_0-\delta,\lambda_0+\delta)$ the equation \eqref{eq4} has a limit cycle $x=\hat x(t;\lambda)$ with $\hat x(t;\lambda_0)=x(t)$. Moreover:
    \begin{itemize}
      \item[\rm{(i.1)}] When $x=x(t)$ is stable (resp. unstable), $x=\hat x(t;\lambda)$ is lower-stable and increases (resp. upper-unstable and decreases) as ${\lambda}$ increases.
      \item[\rm{(i.2)}] When $x=x(t)$ is stable (resp. unstable), $x=\hat x(t;\lambda)$ is upper-stable and decreases (resp. lower-unstable and increases) as ${\lambda}$ decreases.
    \end{itemize}
\item[\rm{(ii)}] If $x=x(t)$ is a semi-stable limit cycle of the equation \eqref{eq4}$|_{\lambda=\lambda_0}$ in the region, then the following statements hold.
\begin{itemize}
  \item [\rm{(ii.1)}] When $x=x(t)$ is upper-unstable and lower-stable (resp. upper-stable and lower-unstable), there exists $\delta>0$ such that for $\lambda\in(\lambda_0-\delta,\lambda_0)$ (resp. $\lambda\in(\lambda_0,\lambda_0+\delta)$) the equation \eqref{eq4} has two limit cycles locating on the distinct sides of $x=x(t)$.
  \item [\rm{(ii.2)}] When $x=x(t)$ is upper-unstable and lower-stable (resp. upper-stable and lower-unstable), the equation \eqref{eq4} has no limit cycle near $x=x(t)$ for $\lambda>\lambda_0$ (resp. $\lambda<\lambda_0$).
\end{itemize}
\end{enumerate}
\end{Proposition}

\section{The maximum multiplicity of limit cycles of equation \eqref{eq2} \label{mA1}}
This section is devoted to giving a criterion for the maximum multiplicity of the non-zero limit cycles of the Abel equation \eqref{eq2} with general coefficients. This will be the last key tool to prove our main result. For the sake of brevity and clarity, in this section we rewrite the equation as
\begin{equation}\label{eq5}
\frac{dx}{dt}=g(t)x^{3}+f(t)x^{2},
\end{equation}
with $g(t)$ and $f(t)$ being the smooth $2\pi$-periodic functions (which can be non-trigonometric).
Since the case when $g(t)$ or $f(t)$ has definite sign is already clear from \cite[Theorem A]{alvarez2007new} (see also section \ref{introduction} and the works \cite{gasull1990limit, lloyd1979note, pliss1966nonlocal}), we mainly analyze the case that both $g(t)$ and $f(t)$ have indefinite signs here.
More precisely, for equation \eqref{eq5} we propose the following two hypotheses:
\begin{enumerate}
\item[{\bf(C.1)}] $f(t)$ and $g(t)$ have no common zeros, and $g(t)$ has exactly two zeros in any $2\pi$-period.
\item[{\bf(C.2)}] The function $u(t):=-\frac{f(t)}{g(t)}$ has the same strict monotonicity in all connected component of $\{t|g(t)\neq 0\}$.
\end{enumerate}
\vskip0.1cm

Now let us consider a non-zero limit cycle $x=x(t)$ of the equation \eqref{eq5}.
For abbreviation we use the notations
\begin{align}\label{eq6}
  x_{*}=\min \limits_{t\in \mathbb R} x(t),\ x^{*}=\max \limits_{t\in \mathbb R} x(t),\ t_{*}\in \{t_{0}|x(t_{0})=x_*\},\ t^{*}\in \{t_{0}|x(t_{0})=x^*\}\cap[t_{*},t_{*}+2\pi),
\end{align}
and also denote by $t_{1}$ and $t_{2}$ the points satisfying
\begin{align}\label{eq7}
 t_1,t_2\in\{t|g(t)=0\}\cap[t_{*},t_{*}+2\pi),\indent t_1<t_2.
\end{align}
 We give the following result.
\begin{Lemma}\label{lem3.1} Suppose that $x=x(t)$ is a non-zero limit cycle of the equation \eqref{eq5}. Then under hypotheses {\bf(C.1)} and {\bf(C.2)}, the following statements hold.
\begin{enumerate}
\item[\rm{(i)}] $x(t)$ has exactly one maximum point, one minimum point and no other stationary point in any $2\pi$-period.
\item[\rm{(ii)}] The point $t^*$ is contained in $(t_1,t_2)$. Moreover, when $u(t)=-\frac{f(t)}{g(t)}$ is strictly monotonically increasing,
$$
\mbox{$x(t)-u(t)$
$\begin{cases}
<0, t\in (t_{*},t_{1})\cup (t^{*},t_{2}),\\
>0, t\in (t_{1},t^{*})\cup (t_{2},t_{*}+2\pi);
\end{cases}$}
$$
and when $u(t)=-\frac{f(t)}{g(t)}$ is strictly monotonically decreasing,
$$
\mbox{$x(t)-u(t)$
$\begin{cases}
>0, t\in (t_{*},t_{1})\cup (t^{*},t_{2}),\\
<0, t\in (t_{1},t^{*})\cup (t_{2},t_{*}+2\pi).
\end{cases}$}
$$
\end{enumerate}
\end{Lemma}
\noindent\bpf \rm{(i)} It is sufficient to show that the conclusion holds in the period $[t_{*},t_{*}+2\pi)$. As defined in \eqref{eq6}, it is clear that $g(t_*)x_*+f(t_*)=0.$ Then, by means of hypothesis {\bf(C.1)}, we have that $t_{*} \in \{t |g(t)\neq 0\}$, and there are $3$ connected components of $g(t)x+f(t)=0$ in the region $[t_{*},t_{*}+2\pi]\times \mathbb R$, that is, $x=u(t)$ with $t \in [t_{*},t_{1}), (t_{1},t_{2})$ and $(t_{2}, t_{*}+2\pi]$.
On the other hand, one can get from hypothesis {\bf(C.2)} that
$$\left | \begin{matrix}
1 & g(t)x^{3}+f(t)x^{2} \\
1 &u^{'}(t)  \\
\end{matrix} \right |_{x=u(t)}  =u^{'}(t)\geq 0\ ( {\rm or}\leq0),$$
which implies that all the orbits of the equation \eqref{eq5} cross the curve $x=u(t)$ from a same side to the other. Thus, the limit cycle $x=x(t)$ can only intersects each connected component of $g(t)x+f(t)=0$ once, and these intersections in the region $[t_{*},t_{*}+2\pi]\times \mathbb R$ must be $(t_{*},x_*)$, $(t^*,x^*)$ and $(t_{*}+2\pi,x_*)$. As a result, $x(t)$ has a unique minimum point $t_*$, a unique maximum point $t^*$ and no other stationary point in $[t_{*},t_{*}+2\pi)$ (more precisely, $t^*\in(t_1,t_2)$). The statement is verified. Also, we additionally obtain that $x(t)$ is increasing (resp. decreasing) in $(t_{*},t^{*})$ (resp. $(t^{*},t_{*}+2\pi)$).
\vskip0.1cm

\rm{(ii)} We only prove the case that $u(t)$ is strictly monotonically increasing, and the other case follows from a similar argument.

According to statement \rm{(i)}, $x=x(t)$ and $x=u(t)$ do not intersect when $t\in(t_{*},t_{*}+2\pi)\backslash \{t_{1},t^{*},t_{2}\}$. Since $x(t_{*})=u(t_{*})$ and $x^{'}(t_{*})=0\leq u^{'}(t_{*})$, we get that $x(t)-u(t)<0$ for $t\in (t_{*},t_{1})$. Similarly, $x(t^{*})=u(t^{*})$ with $x^{'}(t^{*})=0\leq u^{'}(t^{*})$ yields that $x(t)-u(t)>0$ for $t\in (t_{1},t^{*})$ and $x(t)-u(t)<0$ for $t\in (t^{*},t_{2})$. Moreover, note that $x(t_{*}+2\pi)=u(t_{*}+2\pi)$ and $x^{'}(t_{*}+2\pi)=0\leq u^{'}(t_{*}+2\pi)$, we also have $x(t)-u(t)>0$ for $t\in (t_{2},t_{*}+2\pi)$.
This completes the proof.\epf
\vskip0.3cm

In the next we provide our main theorem of this section.
We remark that this result is essentially obtained by an approach established in the work \cite{yu2023number} which studies the equation \eqref{eq5} with some specific $f(t)$ and $g(t)$.

\begin{Theorem}\label{T1}
Under hypotheses {\bf(C.1)} and {\bf(C.2)}, the multiplicity of each non-zero limit cycle of the equation \eqref{eq5} in the region $x\neq0$, is at most two.
Moreover, when the non-hyperbolic limit cycle exists, it must be lower-stable and upper-unstable (resp. lower-unstable and upper-stable) if $u(t)$ is strictly monotonically increasing (resp. decreasing).
\end{Theorem}

\noindent\bpf
We only consider below the case that the function $u(t)$ in Hypothesis {\bf (C.2)} is strictly monotonically increasing, and the opposite case follows exactly
in the same way.

Assume that $x=x(t)$ is a non-zero limit cycle of the equation \eqref{eq5}. Let $x_*$, $x^*$, $t_*$ and $t^*$ be the notations defined as in \eqref{eq6}.
According to Proposition \ref{pro11}, the multiplicity and stability of $x=x(t)$ can be determined by
\begin{equation}\label{eq8}
L'(x_{*})=\exp\big(\int_{t_{*}}^{t_{*}+2\pi} (2f(t)x(t)+3g(t)x^2(t))\, dt\big)
\end{equation}
and
\begin{equation}\label{eq9}
L''(x_{*})=L'(x_{*})\int_{t_{*}}^{t_{*}+2\pi}(2f(t)+6g(t)x(t))\exp\big(\int_{t_{*}}^{t} (2f(s)x(s)+3g(s)x^2(s))\, ds\big)\, dt
\end{equation}
if $L'(x_{*})-1$ and $L''(x_{*})$ do not simultaneously vanish, where $L$ represents the Poincar\'{e} map of the equation \eqref{eq5} from the initial time $t_*$ to the ending time $t_*+2\pi$.

In the following we focus on the case when $x=x(t)$ is non-hyperbolic. First, we shall simply the expression of $L''(x_{*})$. It is clear that $L'(x_*)=1$ in this case. Then,
\begin{equation}\label{eq10}
  \int_{t_{*}}^{t_{*}+2\pi} (2f(t)x(t)+3g(t)x^2(t))\, dt=0,
\end{equation}
and $L''(x_{*})$ is decomposed into $L''(x_{*})=4W_{1}+2W_{2}$, where
\begin{align*}
  &W_{1}=\int_{t_{*}}^{t_{*}+2\pi}g(t)x(t)\exp\big(\int_{t_{*}}^{t} (2f(\tau)x(\tau)+3g(\tau)x^2(\tau))\, d\tau\big)\, dt, \\
  &W_{2}=\int_{t_{*}}^{t_{*}+2\pi}(f(t)+g(t)x(t))\exp\big(\int_{t_{*}}^{t} (2f(\tau)x(\tau)+3g(\tau)x^2(\tau))\, d\tau\big)\, dt.
\end{align*}
Furthermore, note that
\begin{equation}\label{eq11}
 \int_{t_{*}}^{t} \left(f(\tau)x(\tau)+g(\tau)x^2(\tau)\right) d\tau
 =\int_{t_{*}}^{t} \frac{dx(\tau)}{x(\tau)}
 =\ln \frac{x(t)}{x(t_{*})}.
\end{equation}
Setting $h(t):=\int_{t_{*}}^{t} g(\tau)x^2(\tau)\, d\tau$, we get by \eqref{eq10} and \eqref{eq11} that
\begin{align}\label{eq12}
  h(t_{*}+2\pi)
  =\int_{t_{*}}^{t_{*}+2\pi} (2f(t)x(t)+3g(t)x^2(t))\, dt
   -2\int_{t_{*}}^{t_{*}+2\pi} \big(f(t)x(t)+g(t)x^{2}(t)\big)\,dt
  =0.
\end{align}
Therefore, the equalities \eqref{eq11} and \eqref{eq12} yield that
 \begin{align*}
W_{1}&=\int_{t_{*}}^{t_{*}+2\pi} g(t)x(t)\exp\big(2 \ln \frac{x(t)}{x(t_{*})}+h(t)\big)\, dt
    =\frac{1}{x_{*}^{2}}\int_{t_{*}}^{t_{*}+2\pi} g(t)x^{3}(t)\exp h(t)\, dt\\
   & =\frac{1}{x_{*}^{2}}\int_{t_{*}}^{t_{*}+2\pi} x(t) d \exp h(t)
    =-\frac{1}{x_{*}^{2}}\int_{t_{*}}^{t_{*}+2\pi} \exp h(t) d x(t), \\
W_{2}& =\int_{t_{*}}^{t_{*}+2\pi} (f(t)+g(t)x(t))\exp\big(2 \ln \frac{x(t)}{x(t_{*})}+h(t)\big)\, dt\\
   & =\frac{1}{x_{*}^{2}}\int_{t_{*}}^{t_{*}+2\pi} (f(t)x^{2}(t)+g(t)x^{3}(t))\exp h(t)\, dt\\
   & =\frac{1}{x_{*}^{2}}\int_{t_{*}}^{t_{*}+2\pi} \exp h(t) dx(t).
   \label{eq11}
 \end{align*}

Together with statement (i) of Lemma \ref{lem3.1}, the expression of $L''(x_{*})$ can be reduced to
\begin{equation}\label{eq13}
\begin{split}
L''(x_{*})& =-\frac{2}{x_{*}^{2}}\int_{t_{*}}^{t_{*}+2\pi} \exp h(t) d x(t)\\
   & =-\frac{2}{x_{*}^{2}}\big(\int_{t_{*}}^{t^{*}} \exp h(t)d x(t)+\int_{t^{*}}^{t_{*}+2\pi} \exp h(t)d x(t)\big) \\
   & =-\frac{2}{x_{*}^{2}}\big(\int_{x_{*}}^{x^{*}} \exp h(\tau_{1}(x))d x+\int_{x^{*}}^{x_{*}} \exp h(\tau_{2}(x))d x\big) \\
   & =-\frac{2}{x_{*}^{2}}\big(\int_{x_{*}}^{x_{*}} (\exp h(\tau_{1}(x))-\exp h(\tau_{2}(x)))d x\big),
    \end{split}
\end{equation}
where $\tau_{1}$ and $\tau_{2}$ represent the inverse functions of $x|_{[t_{*},t^{*}]}$ and $x|_{[t^{*},t_{*}+2\pi]}$ respectively, satisfying
\begin{equation}\label{eq14}
  \tau_{1}'|_{(x_*,x^*)}>0,\ \tau_{2}'|_{(x_*,x^*)}<0,\ \tau_{1}(x_{*})=t_{*},\ \tau_{2}(x_{*})=t_{*}+2\pi\text{ and } \tau_{1}(x^{*})=\tau_{2}(x^{*})=t^{*}.
\end{equation}

\begin{figure}[!htbp]
\begin{center}
\psfrag{1}{$t_{*}$}
\psfrag{2}{$t_{1}$}
 \psfrag{3}{$0$}
 \psfrag{4}{$t^{*}$}
  \psfrag{5}{$t_{2}$}
  \psfrag{6}{$t_{*}+2\pi$}
   \psfrag{7}{$x_{*}$}
   \psfrag{8}{$x_{1}$}
   \psfrag{9}{$x_{2}$}
   \psfrag{0}{$x^{*}$}
 \psfrag{a}{$u(t)=-f(t)/g(t)$}
  \psfrag{b}{$\sin x=-p_{0}/p_{1}$}
   \psfrag{c}{$x(t)$}
 \psfrag{d}{$t$}
  \psfrag{e}{$x$}
\psfig{file=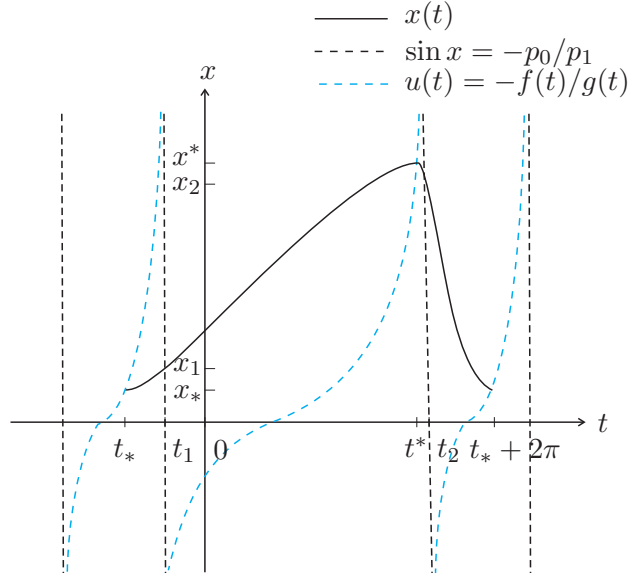,height=3.0in,width=3.3in}
\begin{center}
\caption{Limit cycle $x=x(t)$ of equation \eqref{eq5} when $t\in[t_{*},t_{*}+2\pi)$.}
\label{Fig. 1}
\end{center}
\end{center}
\end{figure}

For convenience of readers, the graph of the limit cycle $x=x(t)$ when it is positive is provided and depicted in Fig. \ref{Fig. 1}.

Next we ascertain the sign of $L''(x_{*})$ by \eqref{eq13}. This can be done once the function
\begin{equation}\label{eq15}
W(s):=h(\tau_{1}(s))-h(\tau_{2}(s))=\int_{\tau_{2}(s)}^{\tau_{1}(s)} g(\tau)x^2(\tau)\, d\tau
\end{equation}
has definite sign in the interval $(x_{*},x^{*})$.
Since, $W(x_*)=W(x^*)=0$ due to \eqref{eq12} and \eqref{eq14}, the problem becomes analyzing the first-order derivative of $W(s)$. We emphasize that $W'(s)$ is continuous in $(x_{*},x^{*})$ with expression
\begin{equation*}
\begin{split}
 W'(s)&=g(\tau_{1}(s))x^{2}(\tau_{1}(s))\tau_{1}^{'}(s)-g(\tau_{2}(s))x^{2}(\tau_{2}(s))\tau_{2}^{'}(s)\\
 &=\frac{g(\tau_{1}(s))}{f(\tau_{1}(s))+g(\tau_{1}(s))x(\tau_{1}(s))}-\frac{g(\tau_{2}(s))}{f(\tau_{2}(s))+g(\tau_{2}(s))x(\tau_{2}(s))} .
\end{split}
\end{equation*}
In addition, let $t_1$ and $t_2$ be defined as in \eqref{eq7}. Observe that $t_{1}\in(t_*,t^{*})$ and $t_{2}\in(t^*,t_*+2\pi)$ from statement (ii) of Lemma \ref{lem3.1}.
Then, $g(\tau_1(s))\neq0$ and $g(\tau_2(s))\neq0$ when $s\not\in\{x(t_1), x(t_2)\}$.
Accordingly, for $s\in (x_{*},x^{*})\backslash \{x(t_{1}),x(t_{2})\}$,
\begin{equation}\label{eq16}
\begin{split}
W'(s) =\hspace{1.5mm}  \frac{1}{x(\tau_{1}(s))-u(\tau_{1}(s))}-\frac{1}{x(\tau_{2}(s))-u(\tau_{2}(s))}.
    \end{split}
\end{equation}

Now denote by $x_{1}=\min \{x(t_{1}),x(t_{2})\}$ and $x_{2}=\max \{x(t_{1}),x(t_{2})\}$. According to the monotonicity of $\tau_1$ and $\tau_2$, we have
\begin{align*}
\left\{
  \begin{aligned}
\tau_{1}(x_{1})\leq t_{1},\\
\tau_{1}(x_{2})\geq t_{1}, \\
\end{aligned}
\right.
\text{ and }
\left\{
  \begin{aligned}
\tau_{2}(x_{1})\geq t_{2},\\
 \tau_{2}(x_{2})\leq t_{2},\\
\end{aligned}
\right.
\end{align*}
respectively.
These yield that
\begin{align*}
  \left\{
  \begin{aligned}
&\tau_{1}(s)\in (t_{*},t_{1}),\ \tau_{2}(s)\in (t_{2},t_{*}+2\pi), &\text{when } s\in (x_{*},x_{1}),\\
&\tau_{1}(s)\in (t_{1},t^{*}),\ \tau_{2}(s)\in (t^{*},t_{2}),      &\text{when } s\in (x_{2}, x^{*}).
  \end{aligned}
  \right.
\end{align*}
Hence, taking statement \rm{(ii)} of Lemma \ref{lem3.1} into account, we get that
\begin{align}\label{eq17A}
x(\tau_1(s))-u(\tau_{1}(s))
\left\{
  \begin{aligned}
<0, s\in (x_{*},x_{1}),\\
>0, s\in (x_{2},x^{*}),
\end{aligned}
  \right.\
\text{ and }\
x(\tau_2(s))-u(\tau_{2}(s))
\left\{
  \begin{aligned}
>0, s\in (x_{*},x_{1}),\\
<0, s\in (x_{2},x^{*}).
\end{aligned}
  \right.
\end{align}

Let us come back to the analysis of the sign of $W^{'}(s)$. The argument is divided into two cases:
\vskip0.2cm

\noindent Case 1. $s\in (x_{*},x^{*})\backslash[x_{1},x_{2}]$.

It follows directly from \eqref{eq16} and \eqref{eq17A} that $W^{'}(s)< 0$ for $s\in (x_{*},x_{1})$ and $W^{'}(s)> 0$ for $s\in (x_{2},x^{*})$.

\vskip0.1cm
\noindent Case 2. $s\in (x_{1},x_{2})$.

 In this case, we write \eqref{eq16} as
$$W'(s) =\frac{V(s)}{\big(x(\tau_{1}(s))-u(\tau_{1}(s))\big)\big(x(\tau_{2}(s))-u(\tau_{2}(s))\big)}$$
with $V(s)=u(\tau_{1}(s))-u(\tau_{2}(s))$. We stress that this expression is well-defined in the interval and therefore $\big(x(\tau_{1}(s))-u(\tau_{1}(s))\big)\big(x(\tau_{2}(s))-u(\tau_{2}(s))\big)\neq0$.
Also, on account of the monotonicity of $\tau_1$, $\tau_2$ and $u$, one can easily see that both the composite functions $u\circ\tau_1$ and $-u\circ\tau_2$ are strictly monotonically increasing in $(x_{1},x_{2})$, and so does $V(s)$. As a result, $W'(s)$ has at most one zero, i.e., it can only change sign at most once, in $(x_{1},x_{2})$.

In summary, according to the continuity of $W'$ and the arguments for the above cases, we confirm that there exists a unique ${x_{0}}\in [x_{1},x_{2}]$, such that $W'({x_{0}})=0$, $W'(x)<0$ for $x\in (x_{*},{x_{0}})$ and $W'(x)>0$ for $x\in ({x_{0}},x^{*})$. This together with $W(x_{*})=W(x^{*})=0$ implies that $W(s)<0$ for $s\in (x_{*},x^{*})$. As a consequence, due to \eqref{eq13}, we finally obtain that $L''(x_{*})>0$. That is to say, the non-zero limit cycle of the equation \eqref{eq5} has multiplicity at most $2$, and it is lower-stable and upper-unstable when it is non-hyperbolic.

The proof of Theorem \ref{T1} is concluded.
\epf

\section{Proof of Theorem \ref{T2}\label{mA2}}
The goal of this section is to prove Theorem \ref{T2} and give a positive answer to the Smale-Pugh-Gasull problem. Let us start from the simplifications for the equation \eqref{eq3}.
By using the transformation $t\mapsto t+\theta$ with $\theta=\arctan\frac{a_2}{a_1}$ (resp. $\theta=\frac{\pi}{2}$) when $a_1\neq 0$ (resp. $a_1=0$),
the equation can be reduced to the following type:
 \begin{equation}\label{eq18}
\frac{dx}{dt}=(p_0+p_1\sin t)x^3+(q_0+q_1\sin t+q_2\cos t)x^2,
\end{equation}
where $p_0=a_0$, $p_1=\sqrt{a_{1}^{2}+a_{2}^{2}}\geq0$, $q_0=b_0$, $q_1=\frac{a_1b_1+a_2b_2}{\sqrt{a_{1}^{2}+a_{2}^{2}}}$ and $q_2=\frac{a_1b_2-a_2b_1}{\sqrt{a_{1}^{2}+a_{2}^{2}}}$.
Moreover, observe that the equation \eqref{eq18} is invariant under the transformations $(t,p_0,q_0,q_2)\mapsto (-t,-p_0,-q_0,-q_2)$ and $(x,t,p_0,q_1)\mapsto (-x,-t,-p_0,-q_1)$, respectively. Thus, without loss of generality, our consideration can be restricted to the case that $q_0, p_0\geq 0$.
On the other hand, thanks to the Theorem \ref{AT1}, it is sufficient to deal with the case when the condition of the third inequality in the theorem is invalid for the equation \eqref{eq18}, that is, $(p_0q_1-p_1q_0)^{2}+p_{0}^{2}q_{2}^{2}<p_{1}^{2}q_{2}^{2}$ (which immediately yields $p_1\neq0$ and $q_2\neq0$).
In summary, we only need to focus on the equation \eqref{eq18} with the following hypothesis:
\vskip0.3cm
\noindent{\bf (H)} $ (p_0q_1-p_1q_0)^{2}+p_{0}^{2}q_{2}^{2}<p_{1}^{2}q_{2}^{2},\ p_0\geq0,\ q_0\geq0,\ p_1>0,\ q_1\in\mathbb R$ and $q_2\neq0$.
\vskip0.3cm

In what follows, we divide our argument into several steps. They are presented as the auxiliary results in the next subsection.
\subsection{Some auxiliary results\label{Some}}
First of all, directly applying Theorem \ref{T1} to the equation \eqref{eq18} yields the following lemma.
\begin{Lemma}\label{lem4.1}
  Under Hypothesis {\bf (H)}, each non-zero limit cycle of the equation \eqref{eq18} has multiplicity at most two. In particular, if such limit cycle is exactly non-hyperbolic, then it must be lower-stable and upper-unstable (resp. lower-unstable and upper-stable) when $q_{2}>0$ (${\rm resp.}\ q_{2}<0$).
\end{Lemma}

\noindent\bpf
  We write the equation \eqref{eq18} as \eqref{eq5} with $g(t)=p_0+p_1\sin t$ and $f(t)=q_0+q_1\sin t+q_2\cos t$. Since the inequality $(p_0q_1-p_1q_0)^{2}+p_{0}^{2}q_{2}^{2}<p_{1}^{2}q_{2}^{2}$ in Hypothesis {\bf (H)} implies that $|p_0|<|p_1|$, the function $g(t)$ has exactly two zeros in any $2\pi$-period.
Furthermore, it follows from a direct calculation that
\begin{align}\label{eq19}
\begin{split}
 \text{sgn}\left(g'(t)f(t)-f'(t)g(t)\right)
 &=\text{sgn}\left(p_{1}q_{2}+p_{0}q_{2}\sin t+(p_{1}q_{0}-p_{0}q_{1})\cos t\right)\\
 &=\text{sgn}\left(p_{1}q_{2}+\sqrt{p_{0}^{2}q_{2}^{2}+(p_{1}q_{0}-p_{0}q_{1})^{2}}\sin(t+\phi)\right)\\
 &=\text{sgn}\left(q_2\right)\\
 &\neq0,
\end{split}
\end{align}
where $\phi=\arctan{\frac{p_{1}q_{0}-p_{0}q_{1}}{p_{0}q_{2}}}$ (resp. $\phi=\frac{\pi}{2}$) when $p_0>0$ (resp. $p_0=0$). Hence, there are no common zeros of $g(t)$ and $f(t)$. The equation \eqref{eq18} satisfies the condition {\bf (C.1)} of Theorem \ref{T1}. Observe that \eqref{eq19} also implies that
\begin{equation*}\label{eq17}
 \text{sgn}\left(\bigg(-\frac{f(t)}{g(t)}\bigg)^{'}\right)
    =\text{sgn}\left(q_2\right).
\end{equation*}
Thus $-\frac{f(t)}{g(t)}$ is strictly monotonically increasing (resp. decreasing) when $q_{2}>0$ (resp. $q_2<0$). The equation also satisfies the condition {\bf (C.2)} of Theorem \ref{T1}. Consequently, our conclusion is obtained on account of Theorem \ref{T1}.
\epf\vskip0.3cm

Next we illustrate the stability of the orbit $x=0$ of the equation \eqref{eq18}. Let $x(t,x_{0})$ be the solution of the equation satisfying the initial condition $x(0,x_{0})=x_{0}$.
It is known that the displacement function $H(x_0):=x(2 \pi ,x_0)-x_0$ is well-defined in a neighbourhood of zero. Then it can be expanded into power series with respect to $x_0$, that is:
\begin{equation*}
  H(x_0) =\sum_{i=2}^{\infty} L_ix_0^i,
\end{equation*}
where the coefficients $L_i$'s are called the {\em Lyapunov constants} of the orbit $x=0$ (see for instance \cite{alvarez2007new}). Following the approach in \cite{lloyd1973number,lloyd1983small} or directly the result in \cite{alvarez2007new}, one can easily get that
 $L_2=2\pi q_{0}$; $L_3=2\pi p_{0}$ when $L_2=0$; and $L_4=p_{1}q_{2}\pi$ when $L_2=L_3=0$.
Hence, under Hypothesis {\bf (H)}, the orbit $x=0$ has multiplicity at most $4$, and its stability with respect to the parameters is summarized in the Table \ref{Tab-zero}.
\begin{table}[!hbp]
 \caption{Multiplicity and stability of the orbit $x=0$ for equation \eqref{eq18}}\label{Tab-zero}
  \centering
 {\footnotesize
 \setlength{\tabcolsep}{5.2mm}
  \begin{tabular}{|c|c|c|}
  \toprule
{Classification of parameters} &Multiplicity  & Stability\\
\hline
  $q_0>0$ &2 &  upper-unstable  and lower-stable\\
\hline
$q_0=0, p_0>0$ &3  &  upper-unstable  and lower-unstable\\
\hline
$p_0=q_0=0, q_2>0$ &4  &upper-unstable  and lower-stable\\
\hline
$p_0=q_0=0, q_2<0$ &4 &upper-stable  and lower-unstable\\
\bottomrule
\end{tabular}}
\end{table}

Now by using Lemma \ref{lem4.1} and the stability of the orbit $x=0$, we provide a preliminary estimate for the number of limit cycles of the equation \eqref{eq18}.
\begin{Lemma}\label{lem4.2}
Assume that Hypothesis {\bf (H)} holds with $q_0\neq0$ for equation \eqref{eq18}.
\begin{itemize}
  \item[\rm{(i)}] When $q_2>0$, equation \eqref{eq18} has at most two limit cycles (taking into account multiplicities) in $x>0$ and $x<0$, respectively.
  \item[\rm{(ii)}] When $q_2<0$, equation \eqref{eq18} has at most one limit cycles (taking into account multiplicities) in $x>0$ and $x<0$, respectively.
\end{itemize}
\end{Lemma}

\noindent\bpf
First we rewrite the equation \eqref{eq18} into the form $\frac{dx}{dt}=S(x, t; q_{0})$ with $S(x, t; q_{0}):=(p_0+p_1\sin t)x^3+(q_0+q_1\sin t+q_2\cos t)x^2$. It is easy to check that $\frac{\partial S}{\partial q_{0}}=x^2>0$ for $x\neq0$. Then from Definition \ref{def1}, the equation \eqref{eq18} forms a rotated equation with respect to the parameter $q_{0}$.

In the following we verify the statement (i) and statement (ii) one by one. We only give the estimates for the number of positive limit cycles of the equation because the conclusions for the negative limit cycles follows exactly from a same argument.
\vskip0.2cm

(i) Assume for a contradiction that there exists $q>0$ such that the equation \eqref{eq18}$|_{q_0=q}$ has $m\geq3$ positive limit cycles (taking into account multiplicities). Without loss of generality, our analysis can be restricted to the fact that all these limit cycles are hyperbolic. Indeed, according to Lemma \ref{lem4.1} and statement (ii.1) of Proposition \ref{pro12}, all the non-hyperbolic limit cycles of the equation split into two hyperbolic limit cycles respectively in a small decrease of $q_0$. Then the equation \eqref{eq18}$|_{q_0=q-\varepsilon}$ can possess only hyperbolic positive limit cycles as $0<\varepsilon\ll q$ and their number is at least $m$, which becomes the case of our concerns.

Now we are able to denote by $x=x_1(t;q_0),\ x=x_2(t;q_0),\ \cdots,\ x=x_m(t;q_0)$ the $m$ positive limit cycles of the equation \eqref{eq18}, with $q_0$ being in some neighborhood of $q$ and $0<x_1(t;q)<\cdots<x_m(t;q)$.
Since $q>0$, we know by Table \ref{Tab-zero} that the orbit $x=0$ is upper-unstable for the equation \eqref{eq18}$|_{q_0=q}$. This yields that both $x=x_{1}(t;q)$ and $x=x_{3}(t;q)$ are stable, and $x=x_{2}(t;q)$ is unstable. Therefore, when $q_{0}$ decreases from $q$, we have by statement (i.2) of Proposition \ref{pro12} that $x=x_{2}(t;q_0)$ is lower-unstable and increases, and $x=x_{3}(t;q_0)$ is upper-stable and decreases, respectively. Furthermore, since Lemma \ref{lem4.1} tells us that the non-hyperbolic limit cycles of \eqref{eq18} must be lower-stable and upper-unstable when $q_2>0$ (including the critical case with $q_0=0$), the limit cycles $x=x_{2}(t;q_0)$ and $x=x_{3}(t;q_0)$ keep their hyperbolicity and stabilities, and do exist for $q_0\in[0,q]$. This means that, the critical case for the equation, \eqref{eq18}$|_{q_0=0}$, has at least two positive limit cycles, which contradicts to Theorem \ref{AT2}. As a result, $m\leq2$ and therefore the equation \eqref{eq18} has at most two positive limit cycles (counting with multiplicities).
\vskip0.2cm

(ii) Similar to the argument in statement (i), it is sufficient to prove that the equation \eqref{eq18} can not possess three consecutive limit cycles $x=0$, $x=x_1(t;q_0)$ and $x=x_2(t;q_0)$ such that the later two are positive and hyperbolic.  In fact, if this case occurs, then it is known by the Table \ref{Tab-zero} that $x=x_1(t;q_0)$ and $x=x_2(t;q_0)$ are stable and unstable, respectively. Thus, on account of statement (i.1) of Proposition \ref{pro12} and Lemma \ref{lem4.1}, the limit cycles $x=x_1(t;q_0)$ and $x=x_2(t;q_0)$ get close to each other with their hyperbolicity and stabilities being retained as $q_0$ increases and satisfies $(p_0q_1-p_1q_0)^{2}+p_{0}^{2}q_{2}^{2}<p_{1}^{2}q_{2}^{2}$. This yields that \eqref{eq18}$|_{(p_0q_1-p_1q_0)^{2}+p_{0}^{2}q_{2}^{2}=p_{1}^{2}q_{2}^{2}}$ has at least two limit cycles (counting with multiplicities), which contradicts to Theorem \ref{AT1}. Accordingly, the conclusion for the number of positive limit cycles of equation \eqref{eq18} is valid.

The proof is finished. \epf\vskip0.3cm

So far, by virtue of Theorem \ref{AT2} and Lemma \ref{lem4.2}, we can actually assert that the equation \eqref{eq18} has at most $4$ non-zero limit cycles under Hypothesis {\bf (H)}. In order to optimize this upper bound and finally prove Theorem \ref{T2}, we provide the last auxiliary result which shows the non-existence of the limit cycles of the equation in some critical cases.

\begin{Lemma}\label{lem4.3}
Assume that $p_0=q_{0}=0$, $p_1>0$ and $q_2\neq0$. Then the equation \eqref{eq18} has no limit cycle in the region $x\neq0$.
\end{Lemma}
\bpf We only prove the case with $q_2>0$ and the rest case with $q_2<0$ follows from a same argument.

By assumption the equation \eqref{eq18} is reduced to
 \begin{equation}\label{eq20}
 \frac{dx}{dt}=S(x,t):=p_1\sin t\cdot x^3+(q_1\sin t+q_2\cos t)\cdot x^2.
\end{equation}
Suppose that $x=x(t)$ is a non-zero limit cycle of the equation \eqref{eq20}. Then we have $\int_{0}^{2 \pi} \big(p_{1}\sin t\cdot x^{2}(t)+(q_{1}\sin t+q_{2}\cos t)\cdot x(t)\big)dt=0$. Therefore,
\begin{equation*}
\begin{split}
\int^{2\pi}_{0}\frac{\partial S}{\partial x}(x(t),t)\ dt
  & =\int_{0}^{2 \pi} \left(3p_{1}\sin t \cdot x^{2}(t)+2(q_{1}\sin t+q_{2}\cos t)\cdot x(t)\right)\ dt\\
  & =p_{1} \int_{0}^{2 \pi} \sin t \cdot x^{2}(t)\,dt\\
  & =p_{1} \left(\int_{0}^{\pi} \sin t \cdot x^{2}(t)\,dt+\int_{\pi}^{2\pi} \sin t \cdot x^{2}(t)\,dt\right)\\
  & =p_{1} \left(\int_{0}^{\pi} \sin t\cdot \left(x^{2}(t)-x^{2}(2\pi-t)\right)dt\right).\\
    \end{split}
\end{equation*}

\begin{figure}[!htbp]
\begin{center}
\psfrag{0}{$0$}
\psfrag{1}{$\pi/2$}
 \psfrag{2}{$\pi$}
  \psfrag{3}{$3\pi/2$}
   \psfrag{4}{$2\pi$}
   \psfrag{5}{$t$}
   \psfrag{6}{$x$}
 \psfrag{7}{$x(2\pi-t)$}
 \psfrag{8}{$x(t)$}
\psfig{file=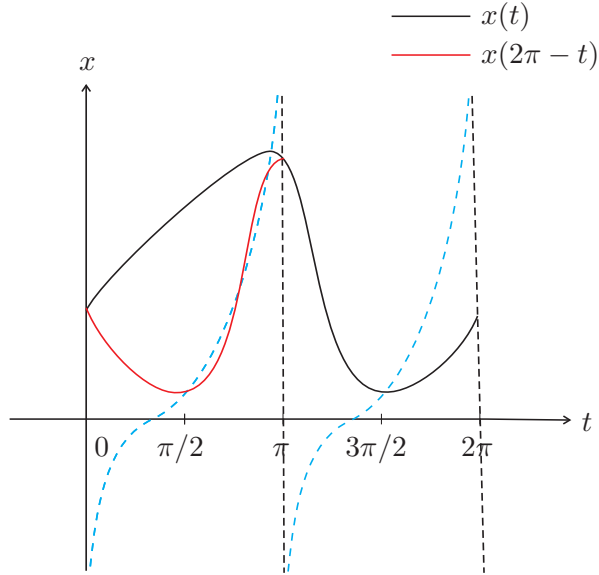,height=3.0in,width=3.2in}
\begin{center}
\caption{The orbit of the equation \eqref{eq20}.}
\label{Fig. 2}
\end{center}
\end{center}
\end{figure}

It can be verified that $v(t):=x(2\pi-t)$ is a solution of the equation
\begin{equation*}
 \frac{dx}{dt}=-S(x,-t)=p_1\sin t\cdot x^3+(q_1\sin t-q_2\cos t)\cdot x^2.
\end{equation*}
Observe that $q_2>0$ implies that $S(x,t)>-S(x,-t)$ (resp. $S(x,t)<-S(x,-t)$) for $t\in[0,\frac{\pi}{2})$ (resp. $t\in(\frac{\pi}{2},\pi]$). Then $x(t)>\nu(t)$ for $t\in(0,\frac{\pi}{2})$ (resp. $t\in(\frac{\pi}{2},\pi)$) taking into account $\nu(0)=x(0)$ (resp. $\nu(\pi)=x(\pi)$), as depicted in Fig. \ref{Fig. 2}.
Consequently, we get that $\int^{2\pi}_{0}\frac{\partial S}{\partial x}(x(t),t)\ dt>0$ (resp. $<0$) when $x(t)>0$ (resp. $<0$).  This is sufficient to show that the limit cycle $x=x(t)$ is unstable (resp. stable) when it is positive (resp. negative).

However, it follows from Table \ref{Tab-zero} that the orbit $x=0$ is upper-unstable and lower-stable when $p_0=q_{0}=0$ and $q_{2}>0$. This contradicts to the stability of $x=x(t)$. As a result, there is no limit cycle in the region $x\neq0$ of the equation \eqref{eq20} and the Theorem is proved (see Fig.~\ref{Fig. 3}(a) and Fig.~\ref{Fig. 4}(a)).
\epf\vskip0.3cm

\subsection{Proof of Theorem \ref{T2}\label{mA2}}
We are ready to prove Theorem \ref{T2}, and provide a positive answer to the Smale-Pugh-Gasull problem now.

\noindent{\bf Proof of Theorem \ref{T2}:}
As has been explained at the beginning of this section, it is sufficient to prove the conclusion for the equation \eqref{eq18} under the Hypothesis {\bf (H)}.
Furthermore, for the subcase $q_2<0$, it is clearly known from Theorem \ref{AT2} (resp. statement (ii) of Lemma \ref{lem4.2}) that the equation \eqref{eq18} has at most three limit cycles (including $x=0$) when $q_0=0$ (resp. $q_0\neq0$). Hence, in the following it only remains to consider the subcase $q_2>0$.

Let $x(t,x_0;p_0,q_0)$ be the solution of the equation \eqref{eq18} with initial value $x_0>0$ and the parameters $p_0$ and $q_0$. According to Lemma \ref{lem4.3} and the stability of $x=0$ from the Table \ref{Tab-zero}, we have that $x(2\pi,x_0;0,0)>x_0$ when $x(2\pi,x_0;0,0)$ is well-defined. In addition, if we denote by $S(x,t;p_0,q_0)$ the right hand side of the equation \eqref{eq18}, then
\begin{align*}
  S(x,t;p_0,q_0)-S(x,t;0,0)=p_0x^3+q_0x^2\geq0,\indent (t,x)\in[0,2\pi]\times\mathbb R^+.
\end{align*}
Thus, $x(2\pi,x_0;p_0,q_0)\geq x(2\pi,x_0;0,0)>x_0$ always holds when $x(2\pi,x_0;p_0,q_0)$ is well-defined, which implies that no positive limit cycle of the equation \eqref{eq18} exists. Consequently, taking Theorem \ref{AT2} and statement (i) of Lemma \ref{lem4.2} into account, the total number of limit cycles of the equation \eqref{eq18}, including $x=0$, is at most three.

\begin{figure}[!htbp]
\centering
\setcounter{subfigure}{0}
\subfigure[$p_{0}=q_{0}=0,q_{2}>0$]
{\psfrag{1}{$t=0$}
 \psfrag{2}{$t=2\pi$}
  \psfrag{3}{$x=0$}
   \psfrag{4}{$x$}
\psfig{file=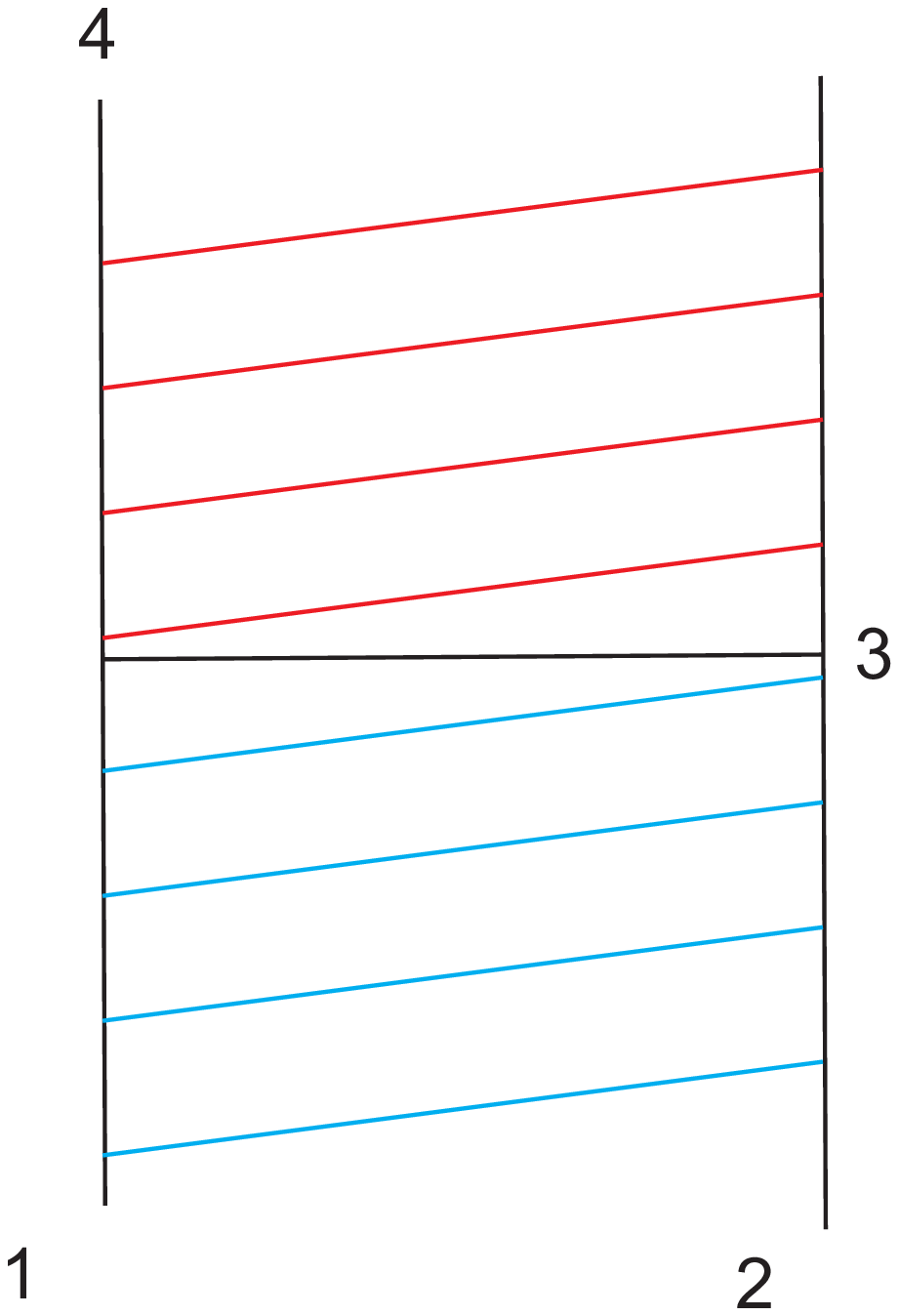,height=1.9in,width=1.5in}}
\hspace{13mm}
\centering
\subfigure[$p_{0}>0,q_{0}=0,q_{2}>0$]
{\psfrag{1}{$t=0$}
 \psfrag{2}{$t=2\pi$}
  \psfrag{3}{$x=0$}
   \psfrag{4}{$x$}
\psfig{file=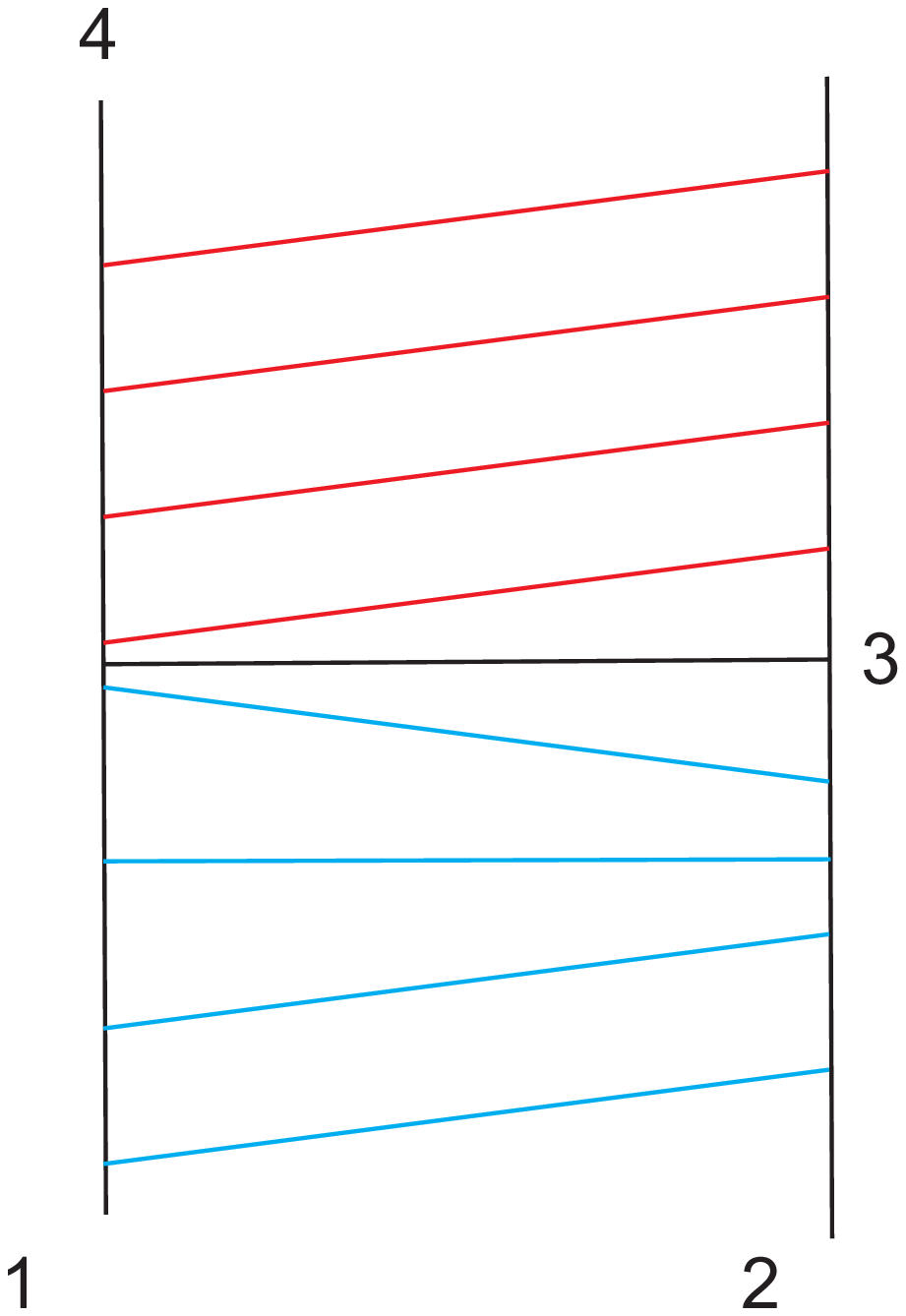,height=1.9in,width=1.5in}}
\hspace{13mm}
\centering
\subfigure[$p_{0},q_{0},q_{2}>0$]
{\psfrag{1}{$t=0$}
 \psfrag{2}{$t=2\pi$}
  \psfrag{3}{$x=0$}
   \psfrag{4}{$x$}
\psfig{file=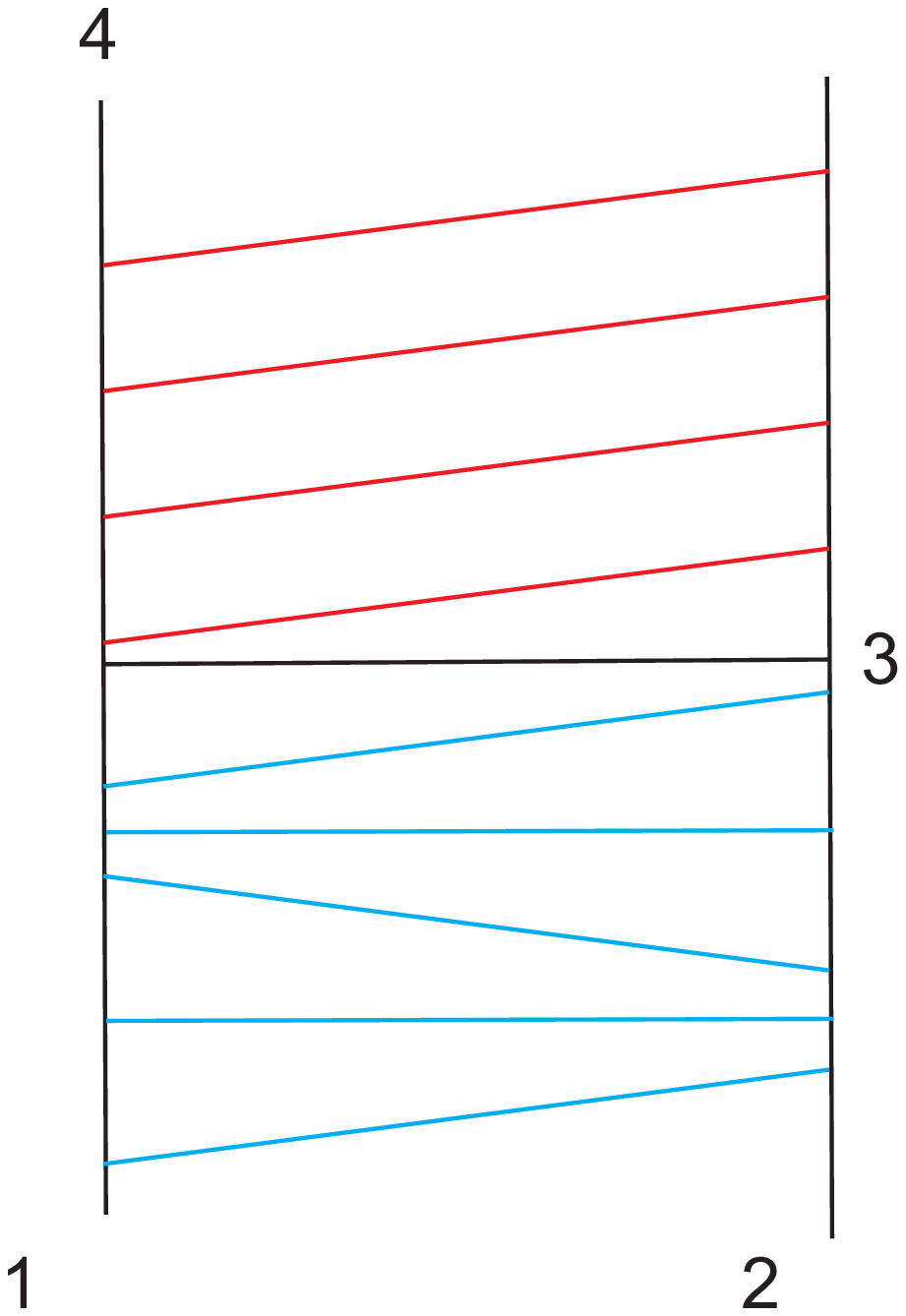,height=1.9in,width=1.5in}}
\begin{center}
\caption{Hopf bifurcations from $x=0$ when $q_{2}>0$.}
\label{Fig. 3}
\end{center}
\end{figure}

\begin{figure}[!htbp]
\centering
\setcounter{subfigure}{0}
\subfigure[$p_{0}=q_{0}=0,q_{2}<0$]
{\psfrag{1}{$t=0$}
 \psfrag{2}{$t=2\pi$}
  \psfrag{3}{$x=0$}
   \psfrag{4}{$x$}
\psfig{file=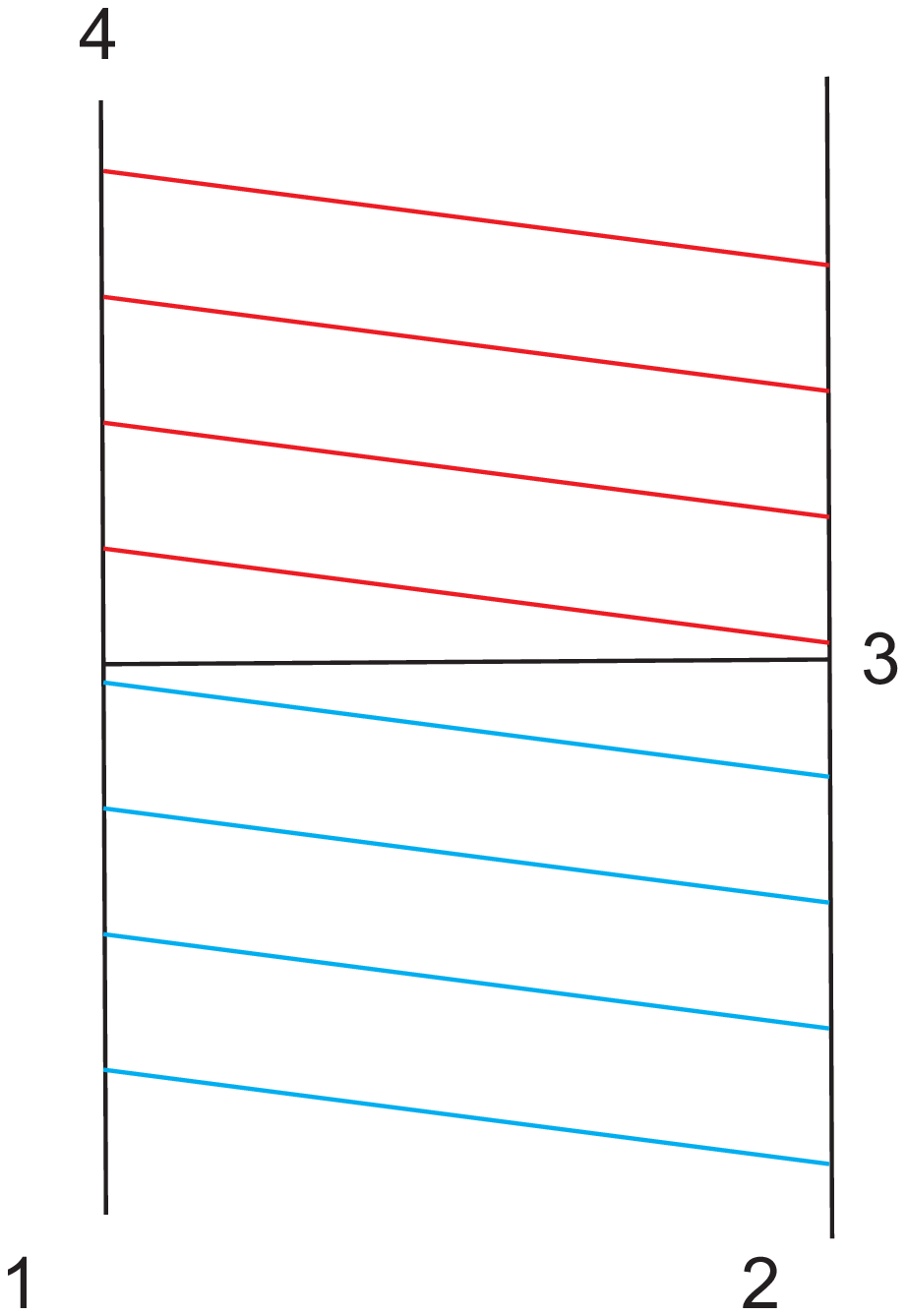,height=1.9in,width=1.5in}}
\hspace{13mm}
\centering
\subfigure[$p_{0}>0,q_{0}=0,q_{2}<0$]
{\psfrag{1}{$t=0$}
 \psfrag{2}{$t=2\pi$}
  \psfrag{3}{$x=0$}
   \psfrag{4}{$x$}
\psfig{file=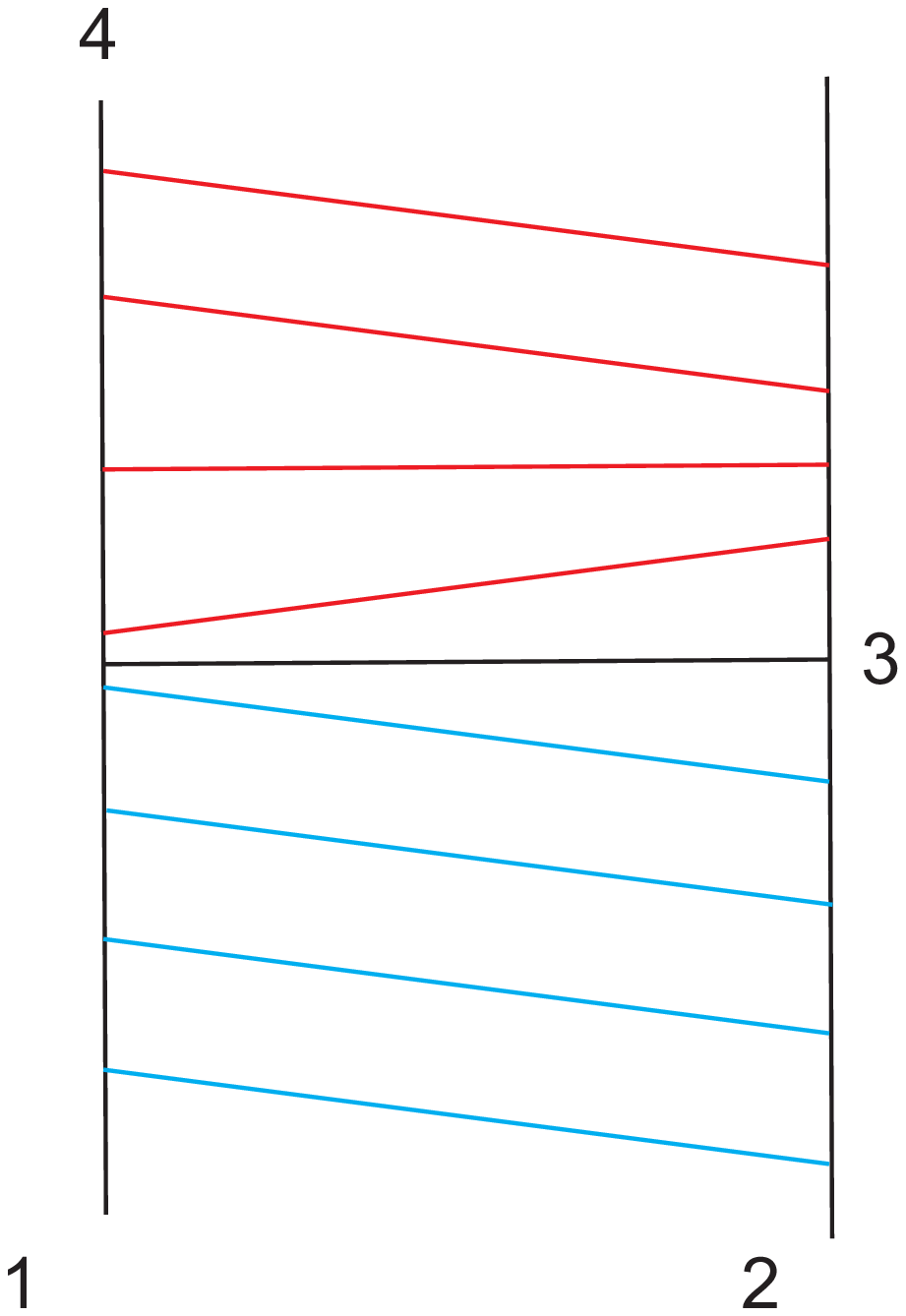,height=1.9in,width=1.5in}}
\hspace{13mm}
\centering
\subfigure[$p_{0},q_{0},q_{2}<0$]
{\psfrag{1}{$t=0$}
 \psfrag{2}{$t=2\pi$}
  \psfrag{3}{$x=0$}
   \psfrag{4}{$x$}
\psfig{file=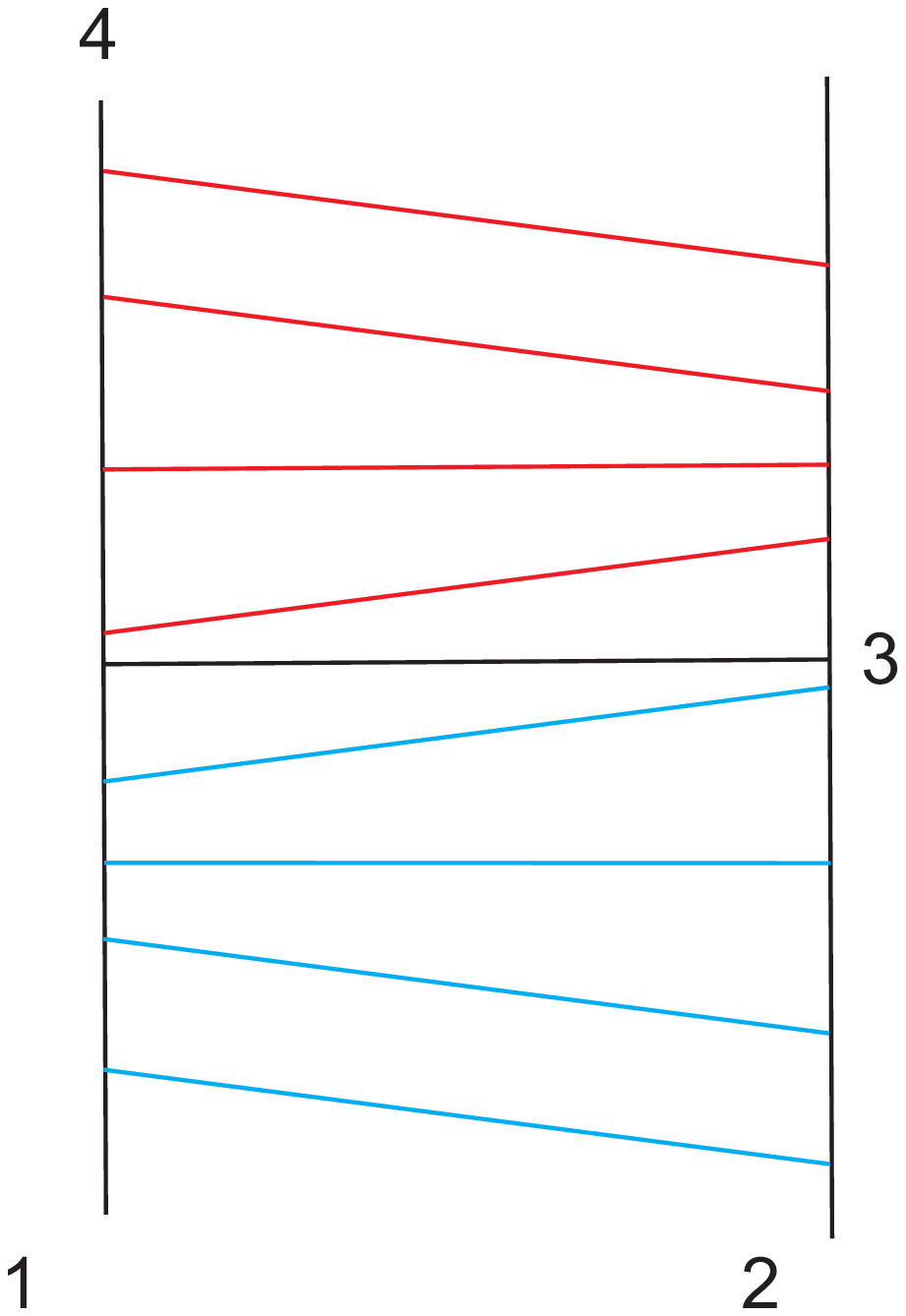,height=1.9in,width=1.5in}}
\begin{center}
\caption{Hopf bifurcations from $x=0$ when $q_{2}<0$.}
\label{Fig. 4}
\end{center}
\end{figure}

Finally we show that the upper bound of the limit cycles can be achieved. It is sufficient to notice again from the Table \ref{Tab-zero} that, the changes of the parameters $p_0$ and $q_0$ yield Hopf bifurcations from $x=0$. The distributions of the limit cycles from these bifurcations are summarized in the following Table \ref{Tab A1}.
For more details see also Fig.~\ref{Fig. 3}(b,c) and Fig.~\ref{Fig. 4}(b,c).

\begin{table}[!hbp]
\caption{Distributions of the limit cycles bifurcating from $x=0$}\label{Tab A1}
  \centering
 {\footnotesize
  \setlength{\tabcolsep}{10mm}
  \begin{tabular}{|c|c|c|}
  \toprule
\backslashbox{Region}{Parameter} & {when $q_2>0$}   & {when $q_2<0$}  \\
\hline
  $x>0$ &0 &1    \\
\hline
$x<0$ &2   &1  \\
\bottomrule
\end{tabular}}
\end{table}

Our proof is finished.
\epf

\section*{Acknowledgements}

The first and third authors are supported by NNSF of China (No. 12171491).
The second author is supported by the NNSF
of China (No. 12271212), and NSF of Guangdong Province (No.2021A1515010029).

\end{document}